# ON INVARIANT MEASURES OF STOCHASTIC RECURSIONS IN A CRITICAL CASE

By Dariusz Buraczewski[1]

*University of Wrocław*

We consider an autoregressive model on $\mathbb{R}$ defined by the recurrence equation $X_n = A_n X_{n-1} + B_n$, where $\{(B_n, A_n)\}$ are i.i.d. random variables valued in $\mathbb{R} \times \mathbb{R}^+$ and $\mathbb{E}[\log A_1] = 0$ (critical case). It was proved by Babillot, Bougerol and Elie that there exists a unique invariant Radon measure of the process $\{X_n\}$. The aim of the paper is to investigate its behavior at infinity. We describe also stationary measures of two other stochastic recursions, including one arising in queuing theory.

**1. Introduction.** We consider the following random process on $\mathbb{R}$:

$$X_n = A_n X_{n-1} + B_n,$$

where the random pairs $(B_n, A_n)$ in $\mathbb{R} \times \mathbb{R}^+$ are independent, identically distributed (i.i.d.) according to a given probability measure $\mu$. This process is called sometimes a first order random coefficients autoregressive model. It appears in various applications, especially in economy and biology; see, for instance [1, 18, 25] and the comprehensive bibliography there.

Many properties of the random process $\{X_n\}$ have been studied under the assumption

(1.1) $$\mathbb{E}[\log A_1] < 0.$$

Then, there exists a unique invariant probability measure $\nu$ of the process $\{X_n\}$ ([17]; see also [12]). This is a measure such that

$$\mu * \nu(f) = \nu(f),$$

---

Received October 2006; revised January 2007.

[1]Supported in part by KBN Grant N201 012 31/1020.

for any positive measurable function $f$ on $\mathbb{R}$, where

$$\mu * \nu(f) = \int_{\mathbb{R} \times \mathbb{R}^+} \int_{\mathbb{R}} f(ax+b)\, d\nu(x)\, d\mu(a,b).$$

Under some further assumptions, the tail behavior of $\nu$ has been described by Kesten [17]. He has proved that

$$(1.2) \qquad \nu(\{t : |t| > x\}) \sim C x^{-\alpha} \qquad \text{as } x \to +\infty$$

for some positive constants $\alpha$ and $C$. Kesten's proof was later essentially simplified by Grincevicius [15] and Goldie [11]; see also [5, 8, 13, 16, 18, 19] for related results.

In this paper we are going to study the so called "critical case," that is, the case when

$$(1.3) \qquad \mathbb{E}[\log A_1] = 0.$$

Then there is no finite invariant measure. However, it has been proved by Babillot, Bougerol and Elie [1] (see also [3, 4]) that there exists a unique (up to a constant factor) stationary Radon measure $\nu$ of the process $\{X_n\}$. Moreover, they have described the behavior of $\nu$ at infinity, proving that, for any given positive $\alpha$ and $\beta$,

$$\nu((\alpha x, \beta x]) \sim \log(\beta/\alpha) \cdot L^{\pm}(|x|) \qquad \text{as } x \to \pm\infty,$$

where $L^{\pm}$ are slowly varying functions.

The aim of this paper is to show that under the hypothesis that the measure $\mu$ is spread out and has some moments, slowly varying functions are constant. Indeed, we show

$$(1.4) \qquad \nu((\alpha x, \beta x]) \sim \log(\beta/\alpha) \cdot C_{\pm} \qquad \text{as } x \to \pm\infty$$

(Theorem 2.2).

If the measure $\mu$ is of a very specific form, that is, $\mu = \mu_t$, where $\{\mu_t\}_{t>0}$ is a one-parameter semi-group of probability measures whose infinitesimal generator is a Hörmander type differential operator on $\mathbb{R}^d \rtimes \mathbb{R}^+$ (or, more generally, on solvable groups of NA type), then the invariant measure $\nu$ has smooth density $m$ and in the situation corresponding to the critical case,

$$(1.5) \qquad x^d m(x \cdot t) \sim c(t) \qquad \text{as } x \to +\infty,$$

for any $t \in \mathbb{R}^d \setminus \{0\}$, where $x \cdot t$ is an appropriate dilation, [6, 7]. Of course, for these particular measures $\mu$, the present result (1.4) follows from (1.5).

In the contractive case (1.1) it was observed by Goldie [11] (see also Grey [13]) that the same problems can be investigated, if the random linear transformations $t \mapsto At + B$ are replaced by a general family of transformations $\{\Psi(t)\}_{t>0}$, provided that, for large values $t$, $\Psi(t)$ is comparable with $At$. Goldie studied several random processes, proved existence of stationary probability measures and described their behavior at infinity that coincided with (1.2). In the critical case we examine a model due to Letac [21]:

$$X'_n = B_n + A_n \max\{X'_{n-1}, C_n\},$$

where $(A_n, B_n, C_n) \in \mathbb{R}^+ \times \mathbb{R} \times \mathbb{R}^+$ are i.i.d. Then, existence and uniqueness of a stationary measure follows from the theory of locally contractive stochastic dynamical systems due to Benda [3]. We shall prove that also in this case the tail behaves regularly and satisfies (1.4) (Theorem 5.1). In particular, the result holds when $B_n = 0$ and we consider the following random process:

$$X''_n = \max\{A_n X''_{n-1}, D_n\},$$

where $(A_n, D_n) \in \mathbb{R}^+ \times \mathbb{R}^+$ are i.i.d. The process $\{X''_n\}$ is called the extremal random process and it plays an important role in modeling of the waiting time for a single server queue [10].

The structure of the paper is as follows. In Section 2 we describe the autoregressive process in the critical case and state our main result, Theorem 2.2. Next, in Section 3 we describe results concerning solutions of the Poisson equation and their asymptotic behavior. In Section 4 we conclude the proof of Theorem 2.2. Finally, in Section 5 we investigate the Letac's model.

## 2. Random difference equation $X_n = A_n X_{n-1} + B_n$.

2.1. *Main theorem.* Given a probability measure $\mu$ on $\mathbb{R} \times \mathbb{R}^+$, consider the following Markov chain on $\mathbb{R}$:

$$\begin{aligned}
X_0 &= 0, \\
X_n &= A_n X_{n-1} + B_n,
\end{aligned} \tag{2.1}$$

where $(B_n, A_n)$ are i.i.d. random variables with values in $\mathbb{R} \times \mathbb{R}^+$, distributed according to $\mu$.

The process $\{X_n\}$ is best defined in the group language. Let $G$ be the "$ax + b$" group, that is, $G = \mathbb{R} \rtimes \mathbb{R}^+$, multiplication being defined by

$$(b, a) \cdot (b', a') = (b + ab', aa').$$

$G$ acts on $\mathbb{R}$ by

$$(b, a) \circ x = ax + b, \qquad (b, a) \in G, \ x \in \mathbb{R}.$$

We sample $(B_n, A_n) \in G$ independently according to a measure $\mu$ and we write

$$X_n = (B_n, A_n) \circ X_{n-1} = (B_n, A_n) \cdots (B_1, A_1) \circ 0.$$

This paper is about the critical case, that is, our main assumption is

$$\mathbb{E} \log A = 0.$$

We are interested in the asymptotic behavior of the (unique) invariant Radon measure of the process $\{X_n\}$, that is, the measure $\nu$ on $\mathbb{R}$ satisfying

$$(2.2) \qquad \mu *_G \nu(f) = \nu(f),$$

for any positive measurable function $f$. Here

$$\mu *_G \nu(f) = \int_{\mathbb{R} \times \mathbb{R}^+} \int_{\mathbb{R}} f(ax + b) \, d\nu(x) \, d\mu(a, b).$$

(We write $*_G$ for the convolution induced by the action of $G$ on $\mathbb{R}$ defined as above, in order to distinguish it with the convolution on $\mathbb{R}$, i.e., denoted by $*_{\mathbb{R}}$.)

Existence and uniqueness of such a measure $\nu$ is due to Babillot, Bougerol and Elie [1] (see also [3, 4] for some comments), who proved the following result:

THEOREM 2.1 ([1]).   *Assume*

$$(2.3) \qquad \mathbb{E} \log A = 0;$$

$$(2.4) \qquad A \not\equiv 1;$$

$$(2.5) \qquad \mathbb{P}[Ax + B = x] < 1 \qquad \text{for all } x \in \mathbb{R};$$

$$(2.6) \qquad \mathbb{E}[(|\log A| + \log^+ |B|)^{2+\delta}] < \infty \qquad \text{for some } \delta > 0.$$

*Then there exists a unique (up to a constant factor) invariant Radon measure $\nu$ on $\mathbb{R}$ of the process $\{X_n\}$. Moreover, if*

$$(2.7) \qquad \begin{aligned} &\text{the closed semigroup generated by} \\ &\text{the support of } \mu \text{ is the whole group } G, \end{aligned}$$

*then there exist two slowly varying functions $L^+$ and $L^-$ on $\mathbb{R}^+$ such that, for any $\alpha, \beta > 0$,*

$$(2.8) \qquad \nu((\alpha x, \beta x]) \sim \log(\beta/\alpha) \cdot L^{\pm}(|x|) \qquad \text{as } x \to \pm\infty.$$

*In particular,*

$$(2.9) \qquad \int_{\mathbb{R}} (1 + |x|)^{-\gamma} \, d\nu(x) < \infty,$$

*for any $\gamma > 0$.*

Our aim is to study more precisely the behavior of $\nu$ at infinity. We shall prove, under additional assumptions, that the functions $L^+$ and $L^-$ are just constants.

Define a probability measure $\mu_A$ on $\mathbb{R}^+$ being the projection of $\mu$ onto the second coordinate $\mu_A = \pi_2(\mu)$, that is, for any Borel set $U$ contained in $\mathbb{R}^+$, we put $\mu_A(U) = \mu(\mathbb{R} \times U)$. Recall that a measure on $\mathbb{R}^+$ is called spread-out if for some $n$, its $n$th convolution power has a nonsingular component relative to the Haar measure on $\mathbb{R}^+$.

Our main result is the following:

THEOREM 2.2.    *Suppose that assumptions* (2.3)–(2.7) *are satisfied and, moreover, there exists a positive constant $\delta$ such that*

(2.10)                              $\mathbb{E}A^\delta < \infty \quad and \quad \mathbb{E}A^{-\delta} < \infty;$

(2.11)                              $\mathbb{E}|B|^\delta < \infty;$

(2.12)                              $\mu_A$ *is spread-out.*

*Then, for any positive numbers $\alpha, \beta$, such that $\alpha < \beta$, we have*

(2.13)
$$\lim_{x \to +\infty} \nu((\alpha x, \beta x]) = \log(\beta/\alpha) \cdot C_+,$$
$$\lim_{x \to -\infty} \nu((\beta x, \alpha x]) = \log(\beta/\alpha) \cdot C_-.$$

2.2. *Sketch of the proof.*   Fix two positive numbers $\alpha < \beta$, and define a function $f$ on $\mathbb{R}$:

(2.14)                         $f(x) = \nu(\alpha e^x, \beta e^x].$

Define a measure $\bar{\mu}$ on $\mathbb{R}$:

$$\bar{\mu}(U) = \mu_A(\{x : -\log x \in U\}),$$

for any Borel set $U$. By (2.3), the mean of $\bar{\mu}$ is equal to 0.

Define convolution of a function $g$ with a measure $\eta$ on $\mathbb{R}$:

$$\eta *_{\mathbb{R}} g(x) = \int_{\mathbb{R}} g(x+y) \eta(dy).$$

The key observation is the following lemma.

LEMMA 2.3.    *The function $f$ satisfies the Poisson equation*

(2.15)                         $\bar{\mu} *_{\mathbb{R}} f(x) = f(x) + \psi(x),$

*where*

(2.16)    $\psi(x) = \displaystyle\int_{\mathbb{R} \times \mathbb{R}^+} \left( \nu\left( \frac{\alpha e^x}{a}; \frac{\beta e^x}{a} \right] - \nu\left( \frac{\alpha e^x - b}{a}; \frac{\beta e^x - b}{a} \right] \right) d\mu(b, a).$

Proof.    We have

$$f(x) = \nu(\alpha e^x, \beta e^x] = \mu *_G \nu(\alpha e^x, \beta e^x]$$

$$= \int_{\mathbb{R} \times \mathbb{R}^+} \int_{\mathbb{R}} \mathbf{1}_{(\alpha e^x, \beta e^x]}(as + b) \, d\nu(s) \, d\mu(b, a)$$

$$= \int_{\mathbb{R} \times \mathbb{R}^+} \nu\left( \frac{\alpha e^x - b}{a}; \frac{\beta e^x - b}{a} \right] d\mu(b, a)$$

and, hence,

$$\bar{\mu} *_{\mathbb{R}} f(x) = f(x) + \bar{\mu} *_{\mathbb{R}} f(x) - \int_{\mathbb{R} \times \mathbb{R}^+} \nu\left( \frac{\alpha e^x - b}{a}; \frac{\beta e^x - b}{a} \right] d\mu(b, a)$$

$$= f(x) + \int_{\mathbb{R} \times \mathbb{R}^+} \left[ \nu\left( \frac{\alpha e^x}{a}; \frac{\beta e^x}{a} \right] - \nu\left( \frac{\alpha e^x - b}{a}; \frac{\beta e^x - b}{a} \right] \right] d\mu(b, a)$$

$$= f(x) + \psi(x),$$

which proves (2.16).   $\square$

The Poisson equation on $\mathbb{R}$ was studied in the 1960s (see next section for more comments). It is well known that if the function $\psi$ is good enough, there exists a formula describing solutions of (2.15). This is not our case, however, in the next section we prove that even if $\psi$ does not satisfy classical hypotheses, under some other assumptions (Theorem 3.3), one can obtain enough information to describe the behavior of the function $f$ at infinity. In Section 4 we prove that $\psi$ possesses all the required properties and conclude the proof of Theorem 2.2.

## 3. Poisson equation.

3.1. *Classical results.*  Given a probability measure $\bar{\mu}$ on $\mathbb{R}$ that is centered, that is, $\int_{\mathbb{R}} x \, d\bar{\mu}(x) = 0$, and a function $\phi$ on $\mathbb{R}$, consider the following Poisson equation:

$$(3.1) \qquad\qquad \bar{\mu} *_{\mathbb{R}} f(x) = f(x) + \phi(x) \qquad dx \text{ a.s.}$$

Positive solutions of this equation were studied in a number of papers (see [2, 22, 24, 26]). For our purpose, we use the result of Port and Stone [24], who have considered the situation when the function $\phi$ is bounded and compactly supported. Then, assuming that the measure $\bar{\mu}$ is spread out, they have found an explicit formula for solutions of (3.1) that are bounded from below. Weaker assumptions on $\bar{\mu}$ were imposed by Ornstein [22] and Baldi [2]; however, they studied only positive functions $\phi$.

Port and Stone [24] define a potential kernel $A$, which can be explicitly computed for $\phi \in \Phi$ ($\Phi$ is the class of bounded measurable, compactly supported functions). The kernel is of the form

$$(3.2) \qquad A\phi(x) = a *_{\mathbb{R}} \phi(x) - \mu_2 *_{\mathbb{R}} \phi(x) + bJ(\phi) - \phi(x), \qquad \phi \in \Phi$$

([24], Theorem 7.1), where $\mu_2$ is a finite positive measure, $b$ is an appropriately chosen constant, $J(\phi) = \int_{\mathbb{R}} \phi(x)\,dx$, $a$ is a continuous function [23] satisfying

$$(3.3) \qquad \lim_{x \to \pm\infty} (a(x-y) - a(x)) = \mp\sigma^{-2}y,$$

where $\sigma^2 = \int_{\mathbb{R}} x^2 \bar{\mu}(dx)$, and moreover, the convergence is uniform w.r.t. $y$ in compact sets.

The potential $A$ provides solutions of the Poisson equation:

**Theorem 3.1** (Port and Stone [24], Theorem 10.3). *Assume that the probability measure $\bar{\mu}$ is centered, spread-out and its second moment is finite. Then, if $\phi$ is an element of $\Phi$, all solutions of Poisson equation* (3.1), *which are bounded from below, are of the form*

$$(3.4) \qquad f(x) = A\phi(x) + \frac{cJ(\phi)}{\sigma^2}x + d, \qquad dx \ a.s.,$$

*where $d$ is any constant and $|c| \leq 1$.*

We observe that (3.2) and (3.3) imply that for large $x$ and $|c| \leq 1$ the absolute value of the term $\frac{cJ(\phi)}{\sigma^2}x$ is dominated by $a *_{\mathbb{R}} \phi$, so $f(x)$ defined by the right-hand side of (3.4) is bounded from below.

The Poisson equation can be solved for a more general class of functions and all solutions are given by a formula that coincides with (3.4). Indeed, the following result holds.

**Theorem 3.2.** *Assume that a function $\phi$ is bounded, continuous and satisfies $\int_{\mathbb{R}} |\phi(x)x|\,dx < \infty$, then $A\phi$ is a well defined continuous function and all solutions of the Poisson equation* (3.1) *bounded from below are of the form* (3.4).

The proof of the foregoing theorem follows the classical path; we present the details of the argument for the reader's convenience in the Appendix.

3.2. *Behavior at infinity.* Now, we are going to present a result that will be the main step in the proof of Theorem 2.2 and will also be used in our study of the asymptotic behavior of the invariant measure $\nu$ for other recursions considered in the paper.

Let $\nu$ be a Radon measure on $\mathbb{R}$. Assume that for every positive $\gamma$ there exists a constant $C$ such that

$$(3.5) \qquad \nu(0,x] \leq C(1+x^{\gamma}), \qquad x \geq 0.$$

Let $f_{\alpha,\beta}(x) = \nu(\alpha e^x, \beta e^x]$. Suppose $f_{\alpha,\beta}$ satisfies the Poisson equation

$$(3.6) \qquad \bar{\mu} *_{\mathbb{R}} f_{\alpha,\beta}(x) = f_{\alpha,\beta}(x) + \psi_{\alpha,\beta}(x), \qquad x \in \mathbb{R},$$

where the measure $\bar{\mu}$ satisfies the assumptions of Theorem 3.1 and also

$$(3.7) \qquad \int_0^\infty e^{\gamma x} \, d\bar{\mu}(x) < \infty$$

for some $\gamma > 0$.

THEOREM 3.3.   *Assume*

$$(3.8) \qquad \int_{-\infty}^\infty |\psi_{\alpha,\beta}(x)x| \, dx < \infty;$$

$$(3.9) \qquad \int_{-1}^1 |\psi_{\alpha,\beta}(x)| \, dx < \infty.$$

*We let*

$$C_{\alpha,\beta}^1 = \int_{\mathbb{R}} \psi_{\alpha,\beta}(x) \, dx,$$

$$C_{\alpha,\beta}^2 = -\int_{\mathbb{R}} (x+1)\psi_{\alpha,\beta}(x) \, dx.$$

*Suppose that*

$$(3.10) \quad \text{the functions } (\alpha,\beta) \mapsto C_{\alpha,\beta}^1 \text{ and } (\alpha,\beta) \mapsto C_{\alpha,\beta}^2 \text{ are continuous.}$$

*If a function $f_{\alpha,\beta}$ satisfies the Poisson equation (3.6), then*

$$\lim_{x \to +\infty} \frac{f_{\alpha,\beta}(x)}{x} = \frac{2C_{\alpha,\beta}^1}{\sigma^2}.$$

*Moreover, if $C_{\alpha,\beta}^1 = 0$, then*

$$\lim_{x \to +\infty} f_{\alpha,\beta}(x) = \frac{2C_{\alpha,\beta}^2}{\sigma^2}.$$

Our aim is to solve the Poisson equation (3.6). Assumptions (3.8)–(3.10) are not sufficient to do it directly; therefore, we are going to use a smoothing operator to obtain a "smoothed" version of (3.6). This will satisfy all the hypotheses of Theorem 3.2. We prove that the solution of the regularized equation behaves regularly at infinity, and finally, we obtain the description of the asymptotic behavior of the function $f_{\alpha,\beta}$.

Consider the function $K(t) = e^{-t}\mathbf{1}_{[0,\infty)}(t)$ and define a smoothing operator (cf. [11, 15]):

$$(3.11) \qquad \breve{g}(t) = K *_{\mathbb{R}} g(t) = \int_{-\infty}^{t} e^{-(t-u)} g(u)\, du.$$

LEMMA 3.4. *The functions $\breve{\psi}_{\alpha,\beta}$, $\breve{f}_{\alpha,\beta}$ and $(\bar{\mu} *_{\mathbb{R}} f_{\alpha,\beta})^{\smile}$ are well-defined. Moreover, $(\bar{\mu} *_{\mathbb{R}} f_{\alpha,\beta})^{\smile} = \bar{\mu} *_{\mathbb{R}} \breve{f}_{\alpha,\beta}$.*

PROOF. By (3.8), $\psi_{\alpha,\beta}$ is integrable, hence, $\breve{\psi}_{\alpha,\beta}$ as the convolution of $K$ with $\psi_{\alpha,\beta}$ is well defined. Next, for every $t$, we have

$$\breve{f}_{\alpha,\beta}(t) = \int_{-\infty}^{t} e^{-(t-u)} \nu(\alpha e^u; \beta e^u]\, du \leq \nu(0; \beta e^t] \int_{-\infty}^{t} e^{-(t-u)}\, du < \infty.$$

Observe, that by (3.5) and the Tonelli theorem, we obtain

$$\begin{aligned}
(\bar{\mu} *_{\mathbb{R}} f_{\alpha,\beta})^{\smile}(t) &= \int_{-\infty}^{t} e^{-(t-u)} \int_{\mathbb{R}} f_{\alpha,\beta}(u+x)\, d\bar{\mu}(x)\, du \\
&= \int_{\mathbb{R}} \int_{-\infty}^{t+x} e^{-(t+x-u)} f_{\alpha,\beta}(u)\, du\, d\bar{\mu}(x) \\
&\leq \int_{\mathbb{R}} \nu(0, \beta e^{t+x}]\, d\bar{\mu}(x) \leq C \int_{\mathbb{R}} (1 + e^{\gamma(t+x)})\, d\bar{\mu}(x),
\end{aligned}$$

which, in view of (3.7), is finite. Finally, the Fubini theorem implies the last conclusion of the lemma. $\square$

Applying the smoothing operator to both sides of (3.6), we obtain a new Poisson equation:

$$(3.12) \qquad \bar{\mu} *_{\mathbb{R}} \breve{f}_{\alpha,\beta}(x) = \breve{f}_{\alpha,\beta}(x) + \breve{\psi}_{\alpha,\beta}(x), \qquad x \in \mathbb{R}.$$

Moreover, in this case the function $\breve{\psi}_{\alpha,\beta}$ is good enough to yield a description of solutions of the equation.

LEMMA 3.5. *The function $\breve{\psi}_{\alpha,\beta}$ satisfies assumptions of Theorem 3.2, and moreover, it vanishes at infinity.*

PROOF. Of course, the function $\breve{\psi}_{\alpha,\beta}$ is continuous. To prove boundedness of $\breve{\psi}_{\alpha,\beta}$, it is enough to prove that its limits at $+\infty$ and $-\infty$ exist. In fact, as we are going to prove, the limits are equal zero. We have

$$\lim_{t \to -\infty} |\breve{\psi}_{\alpha,\beta}(t)| \leq \lim_{t \to -\infty} \int_{-\infty}^{t} |\psi_{\alpha,\beta}(u)|\, du = 0,$$

$$\lim_{t \to +\infty} |\breve{\psi}_{\alpha,\beta}(t)| \leq \lim_{t \to +\infty} \int_{-\infty}^{\infty} e^{-(t-u)} \mathbf{1}_{[0,\infty)}(t-u) |\psi_{\alpha,\beta}(u)| \, du$$

$$= \int_{-\infty}^{\infty} \lim_{t \to +\infty} e^{-(t-u)} \mathbf{1}_{[0,\infty)}(t-u) |\psi_{\alpha,\beta}(u)| \, du = 0.$$

The validity of the equality follows from Lebesgue's theorem. We also have

$$
\begin{aligned}
\int_{-\infty}^{\infty} |\breve{\psi}_{\alpha,\beta}(x) x| \, dx &\leq \int_{-\infty}^{\infty} \int_{-\infty}^{x} |x e^{(t-x)} \psi_{\alpha,\beta}(t)| \, dt \, dx \\
&= \int_{-\infty}^{\infty} e^t |\psi_{\alpha,\beta}(t)| \cdot \left( \int_t^{\infty} |x| e^{-x} \, dx \right) dt \\
&= \int_{-\infty}^{0} e^t |\psi_{\alpha,\beta}(t)| (2 - (t+1) e^{-t}) \, dt \\
&\quad + \int_0^{\infty} (t+1) |\psi_{\alpha,\beta}(t)| \, dt.
\end{aligned}
$$
(3.13)

The value above is finite, by (3.8) and (3.9).  $\square$

PROPOSITION 3.6.  *We have*

$$\lim_{x \to +\infty} \frac{\breve{f}_{\alpha,\beta}(x)}{x} = \frac{2 C_{\alpha,\beta}^1}{\sigma^2}.$$
(3.14)

*Moreover, if $C_{\alpha,\beta}^1 = 0$, then*

$$\lim_{x \to +\infty} \breve{f}_{\alpha,\beta}(x) = \frac{2 C_{\alpha,\beta}^2}{\sigma^2}.$$
(3.15)

PROOF.  The function $\breve{f}_{\alpha,\beta}$ is positive and satisfies the Poisson equation (3.12), moreover, $\breve{\psi}_{\alpha,\beta}$ satisfies the assumptions of Theorem 3.2 (Lemma 3.5), therefore, $A\breve{\psi}_{\alpha,\beta}$ is a well-defined continuous function and

$$\breve{f}_{\alpha,\beta}(x) = A\breve{\psi}_{\alpha,\beta}(x) + \frac{c J(\breve{\psi}_{\alpha,\beta})}{\sigma^2} x + d, \qquad dx \text{ a.e.}$$
(3.16)

Notice that both sides of the foregoing equation are continuous functions, therefore, it is satisfied for all $x \in \mathbb{R}$. Next, observe that because $\nu$ is a Radon measure, for any positive $\varepsilon$, there exists a negative $M$ such that $f_{\alpha,\beta}(u) \leq \varepsilon$ for $u < M$. Therefore, for $x < M$,

$$|\breve{f}_{\alpha,\beta}(x)| \leq \int_{-\infty}^{x} e^{-(x-u)} f_{\alpha,\beta}(u) \, du \leq \varepsilon,$$

which implies

$$\lim_{x \to -\infty} \breve{f}_{\alpha,\beta}(x) = 0.$$

By virtue of (3.3), one can easily prove that

$$(3.17) \qquad \lim_{x \to \pm\infty} \frac{a(x)}{x} = \pm\frac{1}{\sigma^2} \quad \text{and} \quad |a(x-y) - a(x)| \leq C(|y|+1)$$

for some constant $C$ and every $x \in \mathbb{R}$. Therefore, applying (3.16) and Lemma 3.5, we obtain

$$
\begin{aligned}
0 &= \lim_{x \to -\infty} \breve{f}_{\alpha,\beta}(x) = \lim_{x \to -\infty} \left[ A\breve{\psi}_{\alpha,\beta}(x) + \frac{cJ(\breve{\psi}_{\alpha,\beta})}{\sigma^2}x \right] + d \\
&= \lim_{x \to -\infty} \left[ a *_{\mathbb{R}} \breve{\psi}_{\alpha,\beta}(x) - \mu_2 *_{\mathbb{R}} \breve{\psi}_{\alpha,\beta}(x) \right. \\
&\qquad\qquad \left. + bJ(\breve{\psi}_{\alpha,\beta}) - \breve{\psi}_{\alpha,\beta}(x) + \frac{cJ(\breve{\psi}_{\alpha,\beta})}{\sigma^2}x \right] + d \\
&= bJ(\breve{\psi}_{\alpha,\beta}) + d \\
&\qquad + \lim_{x \to -\infty} \left[ \int_{\mathbb{R}} (a(x-y) - a(x))\breve{\psi}_{\alpha,\beta}(y)\,dy + \left( a(x) + \frac{c}{\sigma^2}x \right)J(\breve{\psi}_{\alpha,\beta}) \right] \\
&= bJ(\breve{\psi}_{\alpha,\beta}) + d + \frac{1}{\sigma^2}\int_{\mathbb{R}} \breve{\psi}_{\alpha,\beta}(y)y\,dy + J(\breve{\psi}_{\alpha,\beta}) \cdot \lim_{x \to -\infty} x\left( \frac{a(x)}{x} + \frac{c}{\sigma^2} \right),
\end{aligned}
$$

which implies $J(\breve{\psi}_{\alpha,\beta})(c-1) = 0$ and

$$bJ(\breve{\psi}_{\alpha,\beta}) + d = -\frac{1}{\sigma^2}\int_{\mathbb{R}} \breve{\psi}_{\alpha,\beta}(y)y\,dy.$$

In the same way we compute

$$\lim_{x \to +\infty} \frac{\breve{f}_{\alpha,\beta}(x)}{x} = J(\breve{\psi}_{\alpha,\beta}) \cdot \lim_{x \to +\infty} \left( \frac{a(x)}{x} + \frac{c}{\sigma^2} \right) = \frac{2J(\breve{\psi}_{\alpha,\beta})}{\sigma^2},$$

but one can easily prove $J(\breve{\psi}_{\alpha,\beta}) = J(\psi_{\alpha,\beta}) = C^1_{\alpha,\beta}$, which gives (3.14). If $C^1_{\alpha,\beta} = 0$, then

$$\lim_{x \to +\infty} \breve{f}_{\alpha,\beta}(x) = -\frac{2}{\sigma^2}\int_{\mathbb{R}} \breve{\psi}_{\alpha,\beta}(y)y\,dy,$$

and observing

$$\int_{\mathbb{R}} \breve{\psi}_{\alpha,\beta}(y)y\,dy = \int_{\mathbb{R}} (y+1)\psi_{\alpha,\beta}(y)\,dy,$$

we obtain (3.15), which completes the proof. $\quad\square$

PROOF OF THEOREM 3.3. We have just proved

$$\lim_{x \to +\infty} (x\log x)^{-1}\int_0^x f_{\alpha,\beta}(\log t)\,dt = \frac{2C^1_{\alpha,\beta}}{\sigma^2}.$$

Fix $\varepsilon > 0$ and consider $1 < \delta < 1 + \varepsilon$. We have

$$(\delta - 1) \cdot \frac{f_{\alpha,(1+\varepsilon)\beta}(\log x)}{\log x} \geq (x \log x)^{-1} \int_x^{\delta x} f_{\alpha,\beta}(\log t) \, dt \to \frac{2C^1_{\alpha,\beta}}{\sigma^2} \cdot (\delta - 1).$$

Therefore,

$$(3.18) \qquad \liminf_{x \to +\infty} \frac{f_{\alpha,(1+\varepsilon)\beta}(\log x)}{\log x} \geq \frac{2C^1_{\alpha,\beta}}{\sigma^2}.$$

Analogously, taking $\frac{1}{1+\varepsilon} < \delta < 1$, we prove

$$(3.19) \qquad \limsup_{x \to +\infty} \frac{f_{\alpha,\beta/(1+\varepsilon)}(\log x)}{\log x} \leq \frac{2C^1_{\alpha,\beta}}{\sigma^2}.$$

Comparing (3.18) and (3.19), we obtain

$$\frac{2C^1_{\alpha,\beta/(1+\varepsilon)}}{\sigma^2} \leq \liminf_{x \to +\infty} \frac{f_{\alpha,\beta}(x)}{x} \leq \limsup_{x \to +\infty} \frac{f_{\alpha,\beta}(x)}{x} \leq \frac{2C^1_{\alpha,(1+\varepsilon)\beta}}{\sigma^2}$$

for every $\varepsilon > 0$. Finally, passing in the foregoing inequality with $\varepsilon$ to zero and applying (3.10), we conclude

$$\lim_{x \to +\infty} \frac{f_{\alpha,\beta}(x)}{x} = \frac{2C^1_{\alpha,\beta}}{\sigma^2},$$

which completes the proof. The same argument justifies the second part of the theorem. $\square$

## 4. Proof of Theorem 2.2.

In order to prove Theorem 2.2, it is enough to check that all the assumptions of Theorem 3.3 are satisfied and compute explicitly constants. We begin with a simple lemma.

LEMMA 4.1. *For any $\gamma < 1$, there exists $x_0 \in (0,1)$ such that the measure $\tilde{\nu} = \delta_{x_0} *_{\mathbb{R}} \nu$ satisfies*

$$(4.1) \qquad \int_{-1}^1 \frac{1}{|s|^\gamma} \, d\tilde{\nu}(s) < \infty.$$

PROOF. We have

$$\int_0^1 \int_0^1 \frac{1}{|s - x|^\gamma} \, d\nu(s) \, dx = \int_0^1 \int_0^1 \frac{1}{|s - x|^\gamma} \, dx \, d\nu(s)$$

$$\leq \int_0^1 \int_{-1}^1 \frac{1}{|x|^\gamma} \, dx \, d\nu(s) = C \cdot \nu[0,1] < \infty,$$

therefore, there exists $x_0 \in (0,1)$ such that

$$(4.2) \qquad \int_0^1 \frac{1}{|s - x_0|^\gamma} \, d\nu(s) < \infty,$$

which implies the result.   $\square$

Define the following measure on $G$: $\tilde{\mu} = \delta_{(x_0,1)} *_G \mu *_G \delta_{(-x_0,1)}$, that is, $\tilde{\mu}(U) = \mu((-x_0,1) \cdot U \cdot (x_0,1))$, for Borel sets $U \subset G$. Then $\tilde{\mu}$ satisfies all the assumptions of Theorem 2.2. Hence, there exists a Radon measure on $\mathbb{R}$ invariant under $\tilde{\mu}$ and one can easily check that this is exactly the measure $\tilde{\nu}$ (up to a constant factor) defined in the lemma above. Of course, the behavior of both measures $\nu$ and $\tilde{\nu}$ at infinity is the same. Therefore, it is enough to prove Theorem 2.2 for $\tilde{\nu}$. From now we shall consider the measures $\tilde{\mu}$ and $\tilde{\nu}$ instead of $\mu$ and $\nu$. However, to simplify our notation, we just write $\mu$, $\nu$ and we assume that (4.1) is fulfilled.

LEMMA 4.2.   *For any $\gamma > 0$, there exists a constant $C$ such that*
$$\nu[-x,x] \le C(1+x^\gamma).$$

PROOF.   For any positive $x$, we have
$$\frac{\nu(0,x]}{(1+x)^\gamma} \le \int_0^\infty \frac{1}{(1+y)^\gamma} \, d\nu(y),$$
which, in view of (2.9), is finite.   $\square$

Observe that the function $\psi = \psi_{\alpha,\beta}$ can be written in the form
$$\psi(x) = \psi_\beta(x) - \psi_\alpha(x),$$
where
$$\psi_\beta(x) = \psi_\beta^+(x) - \psi_\beta^-(x)$$
$$= \int_{b \ge 0} \nu\left(\frac{\beta e^x - b}{a}; \frac{\beta e^x}{a}\right] d\mu(b,a) - \int_{b<0} \nu\left(\frac{\beta e^x}{a}; \frac{\beta e^x - b}{a}\right] d\mu(b,a).$$

LEMMA 4.3.   *The function $\psi$ satisfies* (3.8) *and* (3.9).

PROOF.   First we are going to prove
$$(4.3) \qquad \int_0^\infty |\psi_\alpha(x)x| \, dx < \infty.$$

Choose two small positive numbers $\delta_1$ and $\delta_2$ such that $0 < \delta_1 < \delta_2 < \delta$, where $\delta$ is as in the hypothesis of Theorem 2.2 and take $\gamma = \frac{\delta_2 - \delta_1}{2}$. We may assume $\gamma < 1$. We have
$$\int_0^\infty x \psi_\alpha^+(x) \, dx \le C \int_0^\infty e^{\gamma x} \psi_\alpha^+(x) \, dx$$

$$= C \int_{b \geq 0} \int_1^\infty t^{\gamma-1} \int_{(\alpha t - b)/a}^{\alpha t/a} d\nu(s) \, dt \, d\mu(b,a)$$

$$\leq C \int_{b \geq 0} \Big[ \int_{1 \leq as/\alpha} \int_{as/\alpha}^{(as+b)/\alpha} t^{\gamma-1} \, dt \, d\nu(s)$$

$$+ \int_{as/\alpha < 1 \leq (as+b)/\alpha} \int_1^{(as+b)/\alpha} t^{\gamma-1} \, dt \, d\nu(s) \Big] \, d\mu(b,a).$$

Let us denote the integrals above by $I$ and $II$, respectively. We have

$$I \leq C' \int_{b \geq 0} \int_{s \geq 0} ((as+b)^\gamma - (as)^\gamma) \, d\nu(s) \, d\mu(b,a).$$

Observe that for every two positive numbers $p$, $q$ and every number $\varepsilon$ with $0 < \gamma < \varepsilon \leq 1$, there exists a positive constant $C$ such that

$$(p+q)^\gamma - p^\gamma \leq C p^{\gamma - \varepsilon} q^\varepsilon.$$

Applying this inequality, we may dominate the expression above by

$$C \int_{\mathbb{R} \times \mathbb{R}^+} \int_{\mathbb{R}^+} a^{\gamma - \delta_2/2} |s|^{\gamma - \delta_2/2} |b|^{\delta_2/2} \, d\nu(s) \, d\mu(b,a).$$

Finally, by the Schwarz inequality, we obtain

$$I \leq C \int_{\mathbb{R}^+} |s|^{-\delta_1/2} \, d\nu(s) \cdot \left( \int_{\mathbb{R} \times \mathbb{R}^+} a^{-\delta_1} \, d\mu(b,a) \right)^{1/2} \cdot \left( \int_{\mathbb{R} \times \mathbb{R}^+} |b|^{\delta_2} \, d\mu(b,a) \right)^{1/2}$$

and by virtue of (2.10), (2.11), (2.9) and (4.1), the value above is finite. To estimate the second integral, we write

$$II = C \int_{b \geq 0} \int_{as/\alpha < 1 < (as+b)/\alpha} \int_1^{(as+b)/\alpha} t^{\gamma-1} \, dt \, d\nu(s) \, d\mu(b,a)$$

$$= C' \int_{b \geq 0} \int_{(\alpha-b)/a < s < \alpha/a} \left( \left( \frac{as+b}{\alpha} \right)^\gamma - 1 \right) d\nu(s) \, d\mu(b,a)$$

$$\leq C' \int_{b \geq 0} \int_{(\alpha-b)/a < s < \alpha/a} \left( \left( 1 + \frac{b}{\alpha} \right)^\gamma - 1 \right) d\nu(s) \, d\mu(b,a)$$

$$\leq C'' \int_{b \geq 0} b^\gamma \nu\left( \frac{\alpha-b}{a}; \frac{\alpha}{a} \right] d\mu(b,a),$$

the last integral being finite by (2.10), (2.11) and Lemma 4.2. Similar argument can be used to estimate the integral of $\psi_\alpha^-$, which proves (4.3), and moreover, a small modification of the calculations above gives (3.9).

To prove $\int_{-\infty}^0 |x\psi(x)| \, dx < \infty$, we use the fact that $\psi = \bar{\mu} *_{\mathbb{R}} f - f$. By (4.1), we have

$$\int_{-\infty}^0 |x\psi(x)| \, dx = \int_0^1 \frac{|\log t|}{t} \nu(\alpha t; \beta t] \, dt \leq C \int_0^1 \frac{1}{t^{1+\gamma}} \int_{\alpha t \leq s \leq \beta t} d\nu(s) \, dt$$

$$\leq C \int_0^\beta \int_{s/\beta}^{s/\alpha} \frac{1}{t^{1+\gamma}} \, dt \, d\nu(s) \leq C \int_0^\beta s^{-\gamma} \, d\nu(s) < \infty.$$

Similarly,

$$\int_{-\infty}^0 |x\bar{\mu} * f(x)| \, dx \leq C \int_{\mathbb{R}^+} \int_0^1 \frac{1}{t^{1+\gamma}} \nu(\alpha at; \beta at] \, dt \, d\mu_A(a)$$

$$\leq C \int_{\mathbb{R}^+} \int_0^{\beta a} \int_{s/(\beta a)}^{s/(\alpha a)} \frac{1}{t^{1+\gamma}} \, dt \, d\nu(s) \, d\mu_A(a)$$

$$\leq C \int_{\mathbb{R}^+} \int_0^{\beta a} a^\gamma s^{-\gamma} \, ds \, d\nu(s)$$

$$\leq C \int_{\mathbb{R}^+} a^\gamma \Big( \int_0^1 s^{-\gamma} \, d\nu(s) + \nu(1, \beta a] \Big) \, d\mu_A(a),$$

so the integral is finite by (2.10), (4.1) and Lemma 4.2. □

LEMMA 4.4. *We have*

$$C_{\alpha,\beta}^1 = \int_{\mathbb{R}} \psi(x) \, dx = \log(\beta/\alpha) \cdot D_+^1,$$

*where*

$$D_+^1 = -\int_{b \geq 0} \nu\Big( -\frac{b}{a}; 0 \Big] \, d\mu(b,a) + \int_{b < 0} \nu\Big( 0; -\frac{b}{a} \Big] \, d\mu(b,a).$$

PROOF. We write

$$\int_{-\infty}^\infty (\psi_\beta^+(x) - \psi_\alpha^+(x)) \, dx$$

$$= \int_{b \geq 0} \int_0^\infty \frac{1}{t} \Big( \int_{(\beta t - b)/a < s \leq \beta t/a} d\nu(s)$$

$$- \int_{(\alpha t - b)/a < s \leq \alpha t/a} d\nu(s) \Big) \, dt \, d\mu(b,a)$$

$$= \int_{b \geq 0} \Big[ \int_{s > 0} \Big( \int_{as/\beta}^{(as+b)/\beta} - \int_{as/\alpha}^{(as+b)/\alpha} \Big) \frac{1}{t} \, dt \, d\nu(s)$$

$$- \int_{-b/a < s \leq 0} \int_{(as+b)/\beta}^{(as+b)/\alpha} \frac{1}{t} \, dt \, d\nu(s) \Big] \, d\mu(b,a)$$

$$= -\log(\beta/\alpha) \int_{b \geq 0} \nu\Big( -\frac{b}{a}; 0 \Big] \, d\mu(b,a).$$

The second part is established in an analogous way. □

Repeating the foregoing calculations, one can prove the following:

LEMMA 4.5. *We have*

$$-C_{\alpha,\beta}^2 = \int_{-\infty}^{\infty} (x+1)\psi(x)\,dx = \log(\beta/\alpha)(D_+^2 + (1+\log(\alpha\beta))D_+^1),$$

*where*

$$D_+^2 = \left(\int_{b\geq 0}\int_{s>0} + \int_{b<0}\int_{s>-b/a}\right)\log\left(\frac{as}{as+b}\right)d\nu(s)\,d\mu(b,a)$$

$$- \int_{b\geq 0}\int_{-b/a<s\leq 0}\log(as+b)\,d\nu(s)\,d\mu(b,a)$$

$$+ \int_{b<0}\int_{0<s\leq -b/a}\log(as)\,d\nu(s)\,d\mu(b,a).$$

PROOF OF THEOREM 2.2.   In view of Lemmas 4.2, 4.3, 4.4 and 4.5, we apply Theorem 3.3 to prove

$$\lim_{x\to+\infty}\frac{\nu(\alpha e^x, \beta e^x]}{x} = \log(\beta/\alpha)\cdot\frac{2D_+^1}{\sigma^2}.$$

Analogously, we may describe behavior of the measure $\nu$ on the negative half-line. Namely, the function $x\mapsto\nu(-\beta e^x, -\alpha e^x]$ satisfies an appropriate Poisson equation and, reasoning as previously, one can prove that all the hypotheses of Theorem 3.3 are satisfied. Therefore, we obtain

$$\lim_{x\to+\infty}\frac{\nu(-\beta e^x, -\alpha e^x]}{x} = \log(\beta/\alpha)\cdot\frac{2D_-^1}{\sigma^2},$$

where

$$D_-^1 = \int_{b\geq 0}\nu\left(-\frac{b}{a};0\right]d\mu(b,a) - \int_{b<0}\nu\left(0;-\frac{b}{a}\right]d\mu(b,a).$$

Notice $D_-^1 = -D_+^1$. However, these two constants should be nonnegative, so $D_+^1 = D_-^1 = 0$. Hence, applying again Theorem 3.3, we prove (2.13) with

$$C_+ = -\frac{2D_+^2}{\sigma^2}\quad\text{and}\quad C_- = -\frac{2D_-^2}{\sigma^2},$$

for $D_+^2$ defined in Lemma 4.5 and

$$D_-^2 = \left(\int_{b<0}\int_{s<0} + \int_{b\geq 0}\int_{s<-b/a}\right)\log\left(\frac{as}{as+b}\right)d\nu(s)\,d\mu(b,a)$$

$$- \int_{b<0}\int_{0\leq s<-b/a}\log(|as+b|)\,d\nu(s)\,d\mu(b,a)$$

$$+ \int_{b\geq 0}\int_{-b/a\leq s<0}\log(|as|)\,d\nu(s)\,d\mu(b,a). \qquad\square$$

Finally, let us observe that the sum of $C_+$ and $C_-$ can be expressed by a quite simple formula, namely, assuming (4.1), one can prove

$$C_+ + C_- = \frac{2}{\sigma^2} \cdot \int_{\mathbb{R} \times \mathbb{R}^+} \int_{\mathbb{R}} \log \left| \frac{as+b}{as} \right| d\nu(s) \, d\mu(b, a).$$

**5. Model due to Letac and extremal random process.** Goldie [11] has noticed that in the contractive case (i.e., when $\mathbb{E} \log A < 0$) the proof of the asymptotic behavior of the invariant measure for the random difference equation can be written in a much more general setting and that the same ideas can be used for other stochastic recursions (see also Grey [13]). The key observation is that the random transformation $t \to At + B$ can be replaced by some other transformations. Although our arguments use the language of groups, under more restrictive hypotheses, one can obtain analogous results for some other recursions for which the invariant measure is no longer probabilistic.

5.1. *Main theorem.* In this section we are going to consider the following process on $\mathbb{R}$, introduced by Letac [21]:

$$(5.1) \qquad \begin{aligned} X_0 &= 0, \\ X_n &= B_n + A_n \max\{C_n, X_{n-1}\}, \end{aligned}$$

where $(A_n, B_n, C_n)$ are i.i.d. random variables with values in $\mathbb{R}^+ \times \mathbb{R} \times \mathbb{R}^+$, distributed according to a given measure $\mu$. To simplify the notation, we write

$$X_n = \Phi(A_n, B_n, C_n) \circ X_{n-1}.$$

This model has been investigated only in the contractive case, $\mathbb{E} \log A < 0$. Then there exists a unique stationary probability measure $\nu$, that is, a measure satisfying

$$\mu \circ \nu(f) = \nu(f),$$

for any positive measurable function $f$, where

$$\mu \circ \nu(f) = \int_{\mathbb{R}^+ \times \mathbb{R} \times \mathbb{R}^+} \int_{\mathbb{R}} f(\Phi(a, b, c) \circ x) \, d\nu(x) \, d\mu(a, b, c).$$

Under some further assumptions, the tail of $\nu$ was described by Goldie [11].

We are going to study the critical case

$$\mathbb{E}[\log A] = 0.$$

As before, we define the measure $\mu_A$ being the projection of $\mu$ onto the first coordinate: $\mu_A = \pi_A(\mu)$.

Our main result concerning the Markov chain $\{X_n\}$ is the following.

THEOREM 5.1. *Assume that there exists a positive constant $\delta$ such that:*

(5.2)            $\mathbb{E}\log A = 0;$

(5.3)            $A \not\equiv 1;$

(5.4)            $\mathbb{P}[\Phi(A,B,C) \circ x = x] < 1 \qquad \text{for all } x \in \mathbb{R};$

(5.5)            $\mu_A$ *is spread-out;*

(5.6)            $\mathbb{E}A^\delta < \infty \quad \text{and} \quad \mathbb{E}A^{-\delta} < \infty;$

(5.7)            $\mathbb{E}|B|^\delta < \infty \quad \text{and} \quad \mathbb{E}|C|^\delta < \infty;$

(5.8)            $B \geq \delta \qquad a.e.$

*Then there exists a unique (up to a constant factor) invariant measure $\nu$ of the process $\{X_n\}$, and moreover, for every positive $\alpha < \beta$,*

$$\lim_{x \to +\infty} \nu((\alpha x; \beta x]) = \log(\beta/\alpha) \cdot C_+,$$

*where $C_+$ is a positive constant given by the following formula:*

(5.9)      $C_+ = \dfrac{2}{\sigma^2} \displaystyle\int_{\mathbb{R}^+ \times \mathbb{R} \times \mathbb{R}^+} \int_{\mathbb{R}} \log\left(\dfrac{\Phi(a,b,c) \circ s}{as}\right) d\nu(s)\, d\mu(a,b,c),$

*for $\sigma^2 = \mathbb{E}[\log^2 A]$.*

REMARK 5.2. The method we use to prove Theorem 5.1 is quite general and it is not hard to see that also other examples of stochastic recursions can be treated similarly. For example, a special case of the Letac's model, which we call extremal random process, defined by the formula

$$X'_n = \max\{A_n X'_{n-1}, D_n\},$$

where $(A_n, D_n)$ are i.i.d. random variables with values in $\mathbb{R}^+ \times \mathbb{R}^+$, also possesses an invariant Radon measure in the critical case. Under hypotheses analogous to (5.2)–(5.7) and assuming $D \geq \delta$ a.s., one can describe the behavior at infinity of the invariant measure of $\{X'_n\}$.

5.2. *Existence and uniqueness of an invariant measure.* To prove existence and uniqueness of an invariant measure $\nu$ of the process $\{X_n\}$, we apply results of Benda [3], who, using ideas of Babillot, Bougerol and Elie [1], has investigated locally contractive stochastic dynamical systems. Benda has proved that if a stochastic dynamical system $\{Y_n^y\}$ on $\mathbb{R}$ (in fact, he worked in much more general settings) satisfies the following conditions:

- recurrence: for some $y$, the set of accumulation points of $\{Y_n^y(\omega)\}$ is nonempty for almost every trajectory $\omega$;

- contraction: for every compact set $K$ and every couple of starting points $x, y \in \mathbb{R}$,

$$\lim_{n \to +\infty} \mathbf{1}_K(Y_n^y(\omega))|Y_n^y(\omega) - Y_n^x(\omega)| = 0,$$

for almost any trajectory $\omega$,

then there exists a unique (up to a constant) invariant Radon measure $\nu$ of the process $\{Y_n\}$.

PROPOSITION 5.3. *Under the hypotheses of Theorem 5.1, there exists a unique invariant Radon measure $\nu$ of the process $\{X_n\}$.*

PROOF. In view of Benda's result, it is enough to justify that the Markov chain $\{X_n\}$ possesses recurrence and contraction properties.

Define two autoregressive processes

$$\overline{X}_n^x = A_n \overline{X}_{n-1}^x + B_n,$$
$$\overline{X}_0^x = x$$

and

$$\overline{\overline{X}}_n^x = A_n \overline{\overline{X}}_{n-1}^x + B_n + A_n C_n,$$
$$\overline{\overline{X}}_0^x = x.$$

Then both processes satisfy both recurrence and contraction condition [1, 4]. Notice

(5.10) $$\overline{X}_n^x \leq X_n^x \leq \overline{\overline{X}}_n^x,$$

where $X_n^x$ is the process defined as in (5.1), but starts from $x$ instead of 0. The inequality above immediately implies recurrence of $\{X_n\}$. Moreover, for any $x, y \in \mathbb{R}$,

$$|X_n^x - X_n^y| \leq A_1 \cdots A_n |x - y| = |\overline{X}_n^x - \overline{X}_n^y|,$$

hence, for compact sets of the form $K = [0, M]$, for every positive constant $M$ and almost every trajectory $\omega$,

$$|X_n^x(\omega) - X_n^y(\omega)| \cdot \mathbf{1}_K(X_n^x(\omega)) \leq |\overline{X}_n^x(\omega) - \overline{X}_n^y(\omega)| \cdot \mathbf{1}_K(\overline{X}_n^x(\omega)) \to 0,$$

as $n \to \infty$, which yields the contraction property of $\{X_n\}$. $\quad\square$

5.3. *Some properties of $\nu$.*   The main result of this section is the following:

PROPOSITION 5.4.   *For any $\gamma > 0$,*

$$\int_{\mathbb{R}} \frac{1}{1+|x|^{\gamma}}\, d\nu(x) < \infty.$$

In fact, this is the only step in the proof where (5.8) is really needed. In the proof of the analogous result for the random difference equation (2.9), the structure of the group $G$ has been heavily used and the argument cannot be applied in our situation.

For our purpose, we need an explicit formula for $\nu$ (cf. [1]). Define a random walk on $\mathbb{R}$:

(5.11)
$$S_0 = 0,$$
$$S_n = \log(A_1 \cdots A_n), \qquad n \geq 1,$$

and consider the downward ladder times of $S_n$:

(5.12)
$$L_0 = 0,$$
$$L_n = \inf\{k > L_{n-1}; S_k < S_{L_{n-1}}\}.$$

Let $L = L_1$. It is known that

(5.13)
$$-\infty < \mathbb{E}S_L < 0$$

(see Feller [10]). Next, consider the Markov chain $W_n = X_{L_n}$.

LEMMA 5.5.   *There exists a unique invariant probability measure $\nu_L$ of the process $\{W_n\}$.*

PROOF.   Observe that the Markov chain $\{W_n\}$ satisfies the following stochastic recursion:

(5.14)
$$W_1 = X_{L_1},$$
$$W_n = \max\{Z_n, M_n W_{n-1} + Q_n\},$$

where $(M_n, Z_n, Q_n)$ are i.i.d. random variables valued in $\mathbb{R}^+ \times \mathbb{R} \times \mathbb{R}^+$ and, moreover, $M_1 =_d e^{S_L}$, $Z_1 =_d X_L$, $Q_1 =_d \overline{X}_L$, where $=_d$ denotes equality of the corresponding distributions. Therefore, applying (5.13), (5.10) and results of Grincevicius [14] and Elie [9], we obtain

(5.15)
$$-\infty < \mathbb{E}[\log M_1] < 0,$$
$$\mathbb{E}[\log^+ |Z_1|] < \infty,$$
$$\mathbb{E}[\log^+ |Q_1|] < \infty.$$

Notice that to obtain the Letac's model, we may write the recursion (5.14) in a slightly different way:

$$W_n = Q_n + M_n \max\{W_{n-1}, Z'_n\},$$

where

$$Z'_n = \frac{Z_n - Q_n}{M_n}.$$

Then, under assumption (5.15), existence and uniqueness of an invariant probability measure was proved by Letac [21] and Goldie [11]. □

The following lemma can be deduced from [20]. However, we need it in much weaker form than the result proved there and the proof can be considerable simplified, therefore, we give all the details for the reader's convenience.

LEMMA 5.6. *Let $X_i$ be a sequence of i.i.d. real valued random variables such that $\mathbb{E}X_i = 0$. Put $S_n = \sum_{i=1}^{n} X_i$ and define the stopping time $T = \min\{n : S_n > 0\}$, then*

$$\mathbb{E}\left[\sum_{n=1}^{T-1} e^{\gamma S_n}\right] < \infty$$

*for any $\gamma > 0$.*

PROOF. We have

$$\mathbb{E}\left[\sum_{n=1}^{T-1} e^{\gamma S_n}\right] = \sum_{k=0}^{\infty} \mathbb{E}[e^{\gamma S_k}; T > k].$$

Define

$$\phi_\gamma(s) = \sum_{k=0}^{\infty} s^k \mathbb{E}[e^{\gamma S_k}; T > k],$$

then, by Spitzer ([26], p. 181, P5 (a), (c)), we have

$$\phi_\gamma(s) = \frac{1}{1 - \mathbb{E}[s^T e^{\gamma S_{T'}}]} = \sum_{k=0}^{\infty} (\mathbb{E}[s^T e^{\gamma S_{T'}}])^k,$$

where $T' = \min\{n : S_n \geq 0\}$.

Let $\mu'$ be the distribution of $S_{T'}$. Let $Y_1, \ldots, Y_n$ be a sequence of i.i.d. random variables distributed according to $\mu'$. By [10],

$$-\infty < \mathbb{E}Y_i < 0.$$

Define $S'_n = \sum_{i=1}^n Y_i$ and notice

$$\phi_\gamma(1) = \sum_{k=0}^\infty (\mathbb{E}[e^{\gamma S_{T'}}])^k = \sum_{k=0}^\infty \mathbb{E}[e^{\gamma S'_k}] = \int_{-\infty}^0 e^{\gamma x} U^-(dx),$$

where $U^- = \sum (\mu')^{*n}$ is the Green kernel. The renewal theorem [10] implies that, for negative $x$, $U^-((x,0))$ increases linearly, therefore, the expression above is finite.  $\square$

LEMMA 5.7.  *For any positive measurable function $f$,*

$$\nu(f) = C \cdot \int_{\mathbb{R}} \mathbb{E}\left[\sum_{n=0}^{L-1} f(X_n^x)\right] \nu_L(dx).$$

*Moreover,*

$$\operatorname{supp} \nu \subset [\delta, \infty).$$

PROOF.   Denote the right-hand side of the foregoing equation by $\nu_1$. It is enough to prove that $\nu_1$ is a Radon measure and it is invariant of the process $\{X_n\}$. Invariance can be proved as in [1], page 482. Next, observe that, by (5.8), $\Phi(A, B, C) \circ [\delta, \infty) \subset [\delta, \infty)$ a.s., which implies $\operatorname{supp} \nu_L \subset [\delta, \infty)$ and $\operatorname{supp} \nu_1 \subset [\delta, \infty)$. To justify that $\nu_1$ is a Radon measure, notice that, for every positive constant $M$, we have

$$\begin{aligned}
\nu_1(\mathbf{1}_{[0,M]}) &= C \int_{\mathbb{R}} \mathbb{E}\left[\sum_{n=0}^{L-1} \mathbf{1}_{[0,M]}(X_n^x)\right] \nu_L(dx) \\
&\leq C \int_\delta^\infty \mathbb{E}\left[\sum_{n=0}^{L-1} \mathbf{1}_{[0,M]}(e^{S_n} x)\right] \nu_L(dx) \\
&\leq C \mathbb{E}\left[\sum_{n=0}^{L-1} \mathbf{1}_{[0,M/\delta]}(e^{S_n})\right] \\
&\leq \frac{CM}{\delta} \cdot \mathbb{E}\left[\sum_{n=0}^{L-1} e^{-S_n}\right].
\end{aligned}$$

The above value is finite by Lemma 5.6. Finally, because of the uniqueness of $\nu$, we obtain $\nu = C\nu_1$ for some positive constant $C$.  $\square$

PROOF OF PROPOSITION 5.4.   Because of (5.10) and (5.8),

$$\begin{aligned}
(5.16) \quad |X_n^x| &\geq |\overline{X}_n^x| = |A_n \cdots A_1 x + A_n \cdots A_2 B_1 + \cdots + B_n| \\
&\geq A_1 \cdots A_n x, \qquad \text{a.s.}
\end{aligned}$$

for every positive $x$. Therefore, we conclude

$$\nu((1+|x|)^{-\gamma}) = \int_\delta^\infty \mathbb{E}\left[\sum_{n=0}^{L-1} \frac{1}{(1+|X_n^x|)^\gamma}\right] d\nu_L(x)$$

$$\leq C \int_\delta^\infty \mathbb{E}\left[\sum_{n=0}^{L-1} \frac{1}{(A_1 \cdots A_n)^\gamma}\right] d\nu_L(x)$$

and the proposition follows from Lemma 5.6. $\square$

5.4. *Proof of Theorem* 5.1. We shall proceed as in the proof of Theorem 2.2. Fix two positive constants $\alpha$ and $\beta$ and define

$$f(x) = \nu((\alpha e^x, \beta e^x]).$$

Then the function $f$ satisfies the following Poisson equation:

$$\bar{\mu} *_\mathbb{R} f(x) = f(x) + \psi(x),$$

where

$$\psi(x) = \int_{\mathbb{R}^+ \times \mathbb{R} \times \mathbb{R}^+} \int_\mathbb{R} [\mathbf{1}_{(\alpha e^x, \beta e^x]}(as)$$
$$- \mathbf{1}_{(\alpha e^x, \beta e^x]}(b + a \max\{c, s\})] \, d\nu(s) \, d\mu(a, b, c).$$

LEMMA 5.8. *The function $\psi$ satisfies assumptions* (3.8) *and* (3.9).

PROOF. Take $\gamma = \delta/4$, for $\delta$ as in Theorem 5.1. We have

$$\int_0^\infty e^{\gamma x} |\psi(x)| \, dx$$

$$= \int_{\mathbb{R}^+ \times \mathbb{R} \times \mathbb{R}^+} \int_1^\infty \int_\mathbb{R} t^{\gamma-1} |\mathbf{1}_{(\alpha t, \beta t]}(as)$$
$$- \mathbf{1}_{(\alpha t, \beta t]}(b + a \max(c, s))| \, d\nu(s) \, dt \, d\mu(a, b, c)$$

$$= \int_{\mathbb{R}^+ \times \mathbb{R} \times \mathbb{R}^+} \int_{s>c} \int_1^\infty t^{\gamma-1} |\mathbf{1}_{(\alpha t, \beta t]}(as)$$
$$- \mathbf{1}_{(\alpha t, \beta t]}(b + as)| \, d\nu(s) \, dt \, d\mu(a, b, c)$$

$$+ \int_{\mathbb{R}^+ \times \mathbb{R} \times \mathbb{R}^+} \int_{c \geq s} \int_1^\infty t^{\gamma-1} |\mathbf{1}_{(\alpha t, \beta t]}(as)$$
$$- \mathbf{1}_{(\alpha t, \beta t]}(b + ac)| \, d\nu(s) \, dt \, d\mu(a, b, c).$$

The fist integral can be estimated as in the proof of Lemma 4.3. The second one we dominate by

$$\int_{\mathbb{R}^+ \times \mathbb{R} \times \mathbb{R}^+} \int_{c \geq s} \int_1^\infty t^{\gamma-1} (\mathbf{1}_{(\alpha t, \beta t]}(b + ac) + \mathbf{1}_{(\alpha t, \beta t]}(as)) \, dt \, d\nu(s) \, d\mu(a, b, c)$$

$$\leq \int_{\mathbb{R}^+ \times \mathbb{R} \times \mathbb{R}^+} \int_{c \geq s} \Big( \int_{(b+ac)/\beta \leq t < (b+ac)/\alpha} t^{\gamma-1} \, dt$$

$$+ \int_{as/\beta \leq t < as/\alpha} t^{\gamma-1} \, dt \Big) \, d\nu(s) \, d\mu(a,b,c)$$

$$\leq C \int_{\mathbb{R}^+ \times \mathbb{R} \times \mathbb{R}^+} \int_{c \geq s} ((b+ac)^\gamma + (as)^\gamma) \, d\nu(s) \, d\mu(a,b,c)$$

$$\leq C \int_{\mathbb{R}^+ \times \mathbb{R} \times \mathbb{R}^+} \int \Big( (b+ac)^\gamma \Big( \frac{c}{s} \Big)^\gamma + (ac)^\gamma \Big( \frac{c}{s} \Big)^\gamma \Big) \, d\nu(s) \, d\mu(a,b,c)$$

$$\leq C \int_{\mathbb{R}^+} s^{-\gamma} \, d\nu(s) \cdot \int_{\mathbb{R}^+ \times \mathbb{R} \times \mathbb{R}^+} (|b|^\gamma c^\gamma + a^\gamma c^{2\gamma}) \, d\mu(a,b,c).$$

Finally, we use the Schwarz inequality and conclude finiteness of the integral. The remaining part of the Lemma can be proved as in lemma 4.3.  $\square$

LEMMA 5.9.   *We have*

$$\int_{\mathbb{R}} \psi(x) \, dx = 0.$$

PROOF.   Applying the Fubini theorem, we have

$$\int_{\mathbb{R}} \psi(x) \, dx$$

$$= \int_{b > \delta} \int_{s > \delta} \Big[ \int \frac{1}{t} \cdot \mathbf{1}_{(\alpha t, \beta t]}(as) \, dt$$

$$- \int \frac{1}{t} \cdot \mathbf{1}_{(\alpha t, \beta t]}(\Phi(a,b,c) \circ s) \, dt \Big] \, d\nu(s) \, d\mu(a,b,c)$$

$$= \int_{b > \delta} \int_{s > \delta} \Big[ \int_{as/\beta \leq t < as/\alpha} \frac{1}{t} \, dt$$

$$- \int_{(\Phi(a,b,c) \circ s)/\beta \leq t < (\Phi(a,b,c) \circ s)/\alpha} \frac{1}{t} \, dt \Big] \, d\nu(s) \, d\mu(a,b,c) = 0.$$

$\square$

LEMMA 5.10.   *We have*

$$\int x\psi(x) \, dx = \log(\beta/\alpha) \cdot \int \log \Big( \frac{as}{\Phi(a,b,c) \circ s} \Big) \, d\nu(s) \, d\mu(a,b,c).$$

PROOF OF THEOREM 5.1.   In view of previous lemmas, the result follows immediately from Theorem 3.3.  $\square$

## APPENDIX: PROOF OF THEOREM 3.2

LEMMA A.1. *The potential $A$, given by formula (3.2), can be defined for every function $\psi$ satisfying assumptions of Theorem 3.2. Moreover, the function $A\psi$ is continuous.*

PROOF. Notice first that $\mu_2 * \psi$ is a continuous function, because $\psi$ is bounded and continuous. Next, $a * \psi(x)$ is finite for every $x$ because, by (3.17),

$$
\begin{aligned}
|a *_{\mathbb{R}} \psi(x)| &\leq \int_{\mathbb{R}} |\psi(y)||a(x-y)|\,dy \\
&\leq C\left(\int_{\mathbb{R}} |\psi(y)|\,dy + |x|\int_{\mathbb{R}} |\psi(y)|\,dy + \int_{\mathbb{R}} |\psi(y)y|\,dy\right) \\
&\leq C(\psi)(1+|x|) < \infty.
\end{aligned}
\tag{A.1}
$$

Finally, to prove continuity of the function $a * \psi$, fix $x \in \mathbb{R}$ and consider a sequence $\{x_n\}$ tending to $x$. Put $h_n(y) = \psi(y)a(x_n - y)$. Then by (3.17), all the functions $|h_n|$ are dominated by $h(y) = |\psi(y)| \cdot (|a(x-y)| + C)$, for an appropriate large constant $C$, which is an integrable function. Therefore, by the Lebesgue theorem and using continuity of the function $a$,

$$
\begin{aligned}
\lim_{n \to +\infty} a *_{\mathbb{R}} \psi(x_n) &= \lim_{n \to +\infty} \int_{\mathbb{R}} h_n(y)\,dy \\
&= \int_{\mathbb{R}} \psi(y)a(x-y)\,dy \\
&= a *_{\mathbb{R}} \psi(x). \qquad \square
\end{aligned}
$$

LEMMA A.11. *For every function $\psi$ satisfying assumptions of Theorem 3.2 and $x \in \mathbb{R}$, the following Poisson equation is fulfilled:*

$$
\bar{\mu} *_{\mathbb{R}} A\psi(x) = A\psi(x) + \psi(x).
$$

PROOF. The foregoing Poisson equation is satisfied when $\psi$ is an element of $\Phi$ (Port and Stone [24], Theorem 10.1). Without any loss of generality, we may assume $\psi \geq 0$. Then take any sequence $\psi_n$ of positive continuous, compactly supported functions, tending pointwise to $\psi$ and satisfying $\psi_n(x) \leq \psi(x)$ for every $x$ and $n$. Then

$$
\bar{\mu} *_{\mathbb{R}} A\psi_n(x) = A\psi_n(x) + \psi_n(x).
$$

Notice that $A\psi_n$ tends pointwise to $A\psi$, which is just a consequence of the Lebesgue theorem. Finally, observe, that by (A.1), $A\psi(x)$ can be bounded by $C(|x|+1)$ for some constant $C$, therefore, $\bar{\mu} * A\psi$ is well defined and applying again the Lebesgue theorem, we conclude that $\bar{\mu} * A\psi_n$ tends to $\bar{\mu} * A\psi$. $\square$

Let us fix a function $f$, which is a solution of the Poisson equation (3.1), for some function $\psi$ satisfying hypotheses of Theorem 3.2. Suppose $g$ is a continuous, compactly supported function, such that $J(g) = 1$. Define a function

$$(A.2) \qquad h(x) = f(x) + J(\psi)Ag(x) - A\psi(x), \qquad x \in \mathbb{R}.$$

LEMMA A.12.  *If the function $f$ is bounded from below, then also the function $h$ is bounded from below.*

PROOF.  It is enough to show that $J(\psi)Ag - A\psi$ is a bounded function. By Lemma A.1, it is a continuous function, therefore, it is enough to justify that the limits

$$\lim_{x \to \pm\infty} (a *_{\mathbb{R}} (J(\psi)g - \psi)(x))$$

exist and are finite. We have, by (3.17),

$$\lim_{x \to +\infty} (a *_{\mathbb{R}} (J(\psi)g - \psi)(x))$$

$$= \lim_{x \to +\infty} \int_{\mathbb{R}} (a(x-y) - a(x))(J(\psi)g(y) - \psi(y)) \, dy.$$

Observe that, because of (3.17),

$$|a(x-y) - a(x)||J(\psi)g(y) - \psi(y)|$$

$$\leq C(|y| + 1)|J(\psi)g(y) - \psi(y)|,$$

which is an integrable function. Therefore, by the Lebesgue theorem, we obtain

$$\lim_{x \to +\infty} (a *_{\mathbb{R}} (J(\psi)g - \psi)(x))$$

$$= -\sigma^{-2} \int_{\mathbb{R}} y(J(\psi)g(y) - \psi(y)) \, dy < \infty,$$

which proves the lemma. $\square$

PROOF OF THEOREM 3.2.  Define $h$ as in (A.2), then $h$ is bounded from below (Lemma A.12) and by Lemma A.11 satisfies

$$\bar{\mu} *_{\mathbb{R}} h(x) = h(x) + J(\psi)g(x), \qquad x \in \mathbb{R}.$$

Therefore, by Theorem 3.1,

$$h(x) = J(\psi)Ag(x) + \frac{cJ(\psi)}{\sigma^2}x + d, \qquad dx \text{ a.s.}$$

Hence,

$$(A.3) \quad \begin{aligned} f(x) &= h(x) - J(\psi)Ag(x) + A\psi(x) \\ &= A\psi(x) + \frac{cJ(\psi)}{\sigma^2} + d, \qquad dx \text{ a.s.} \end{aligned} \qquad \square$$

**Acknowledgments.** This paper was prepared when the author was staying at Department of Mathematics, Université de Rennes and at Department of Mathematics, University Pierre & Marie Curie, Paris VI. The visits were financed by the *European Commission* IHP Network 2002–2006 *Harmonic Analysis and Related Problems* (Contract Number: HPRN-CT-2001-00273 - HARP) and European Commission Marie Curie Host Fellowship for the Transfer of Knowledge "Harmonic Analysis, Nonlinear Analysis and Probability," MTKD-CT-2004-013389. The author would like to express his gratitude to the hosts for hospitality.

The author is grateful to Philippe Bougerol, Sara Brofferio, Ewa Damek, Andrzej Hulanicki and Yves Guivarc'h for valuable discussions on the subject of the paper.

## REFERENCES

[1] Babillot, M., Bougerol, P. and Elie, L. (1997). The random difference equation $X_n = A_n X_{n-1} + B_n$ in the critical case. *Ann. Probab.* **25** 478–493. MR1428518

[2] Baldi, P. (1974). Sur l'équation de Poisson. *Ann. Inst. H. Poincaré Sect. B (N.S.)* **10** 423–434. MR0383550

[3] Benda, M. (1999). Contractive stochastic dynamical systems. Preprint.

[4] Brofferio, S. (2003). How a centred random walk on the affine group goes to infinity. *Ann. Inst. H. Poincaré Probab. Statist.* **39** 371–384. MR1978985

[5] Buraczewski, D., Damek, E., Guivarc'h, Y., Hulanicki, A. and Urban, R. (2007). Tail-homogeneity of stationary measures for some multidimensional stochastic recursions. Unpublished manuscript.

[6] Buraczewski, D., Damek, E. and Hulanicki, A. (2006). Asymptotic behavior of Poisson kernels NA groups. *Comm. Partial Differential Equations* **31** 1547–1589. MR2273965

[7] Damek, E. and Hulanicki, A. (2006). Asymptotic behavior of the invariant measure for a diffusion related to an $NA$ group. *Colloq. Math.* **104** 285–309. MR2197079

[8] de Saporta, B., Guivarc'h, Y. and Le Page, E. (2004). On the multidimensional stochastic equation $Y_{n+1} = A_n Y_n + B_n$. *C. R. Math. Acad. Sci. Paris* **339** 499–502. MR2099549

[9] Élie, L. (1982). Comportement asymptotique du noyau potentiel sur les groupes de Lie. *Ann. Sci. École Norm. Sup. (4)* **15** 257–364. MR0683637

[10] Feller, W. (1966). *An Introduction to Probability Theory and Its Applications* II. Wiley, New York. MR0210154

[11] Goldie, C. M. (1991). Implicit renewal theory and tails of solutions of random equations. *Ann. Appl. Probab.* **1** 126–166. MR1097468

[12] Goldie, C. M. and Maller, R. A. (2000). Stability of perpetuities. *Ann. Probab.* **28** 1195–1218. MR1797309

[13] GREY, D. R. (1994). Regular variation in the tail behaviour of solutions of random difference equations. *Ann. Appl. Probab.* **4** 169–183. MR1258178

[14] GRINCEVIČJUS, A. K. (1975). Limit theorem for products of random linear transformations of the line. *Litovsk. Mat. Sb.* **15** 61–77, 241. MR0413216

[15] GRINCEVIČJUS, A. K. (1975). On a limit distribution for a random walk on lines. *Litovsk. Mat. Sb.* **15** 79–91, 243. MR0448571

[16] GUIVARC'H, Y. (2006). Heavy tail properties of multidimensional stochastic recursions. In *Dynamics and Stochastics*: *Festschrift in Honor of M.S. Keane* (D. Denteneer, F. Den Hollander and E. Verbitskiy, eds.) 85–99. IMS, Beachwood, OH.

[17] KESTEN, H. (1973). Random difference equations and renewal theory for products of random matrices. *Acta Math.* **131** 207–248. MR0440724

[18] KLÜPPELBERG, C. and PERGAMENCHTCHIKOV, S. (2004). The tail of the stationary distribution of a random coefficient AR($q$) model. *Ann. Appl. Probab.* **14** 971–1005. MR2052910

[19] LE PAGE, É. (1983). Théorèmes de renouvellement pour les produits de matrices aléatoires. Équations aux différences aléatoires. In *Séminaires de Probabilités Rennes 1983* 1–116. Univ. Rennes I, Rennes. MR0863321

[20] LE PAGE, É. and PEIGNÉ, M. (1997). A local limit theorem on the semi-direct product of $\mathbf{R}^{*+}$ and $\mathbf{R}^d$. *Ann. Inst. H. Poincaré Probab. Statist.* **33** 223–252. MR1443957

[21] LETAC, G. (1986). A contraction principle for certain Markov chains and its applications. In *Random Matrices and Their Applications* (J. E. Cohen, H. Kesten and C. M. Newman, eds.) 263–273. Amer. Math. Soc., Providence, RI. MR0841098

[22] ORNSTEIN, D. S. (1969). Random walks. I, II. *Trans. Amer. Math. Soc.* **138** 1–43; ibid. **138** 45–60. MR0238399

[23] PORT, S. C. and STONE, C. J. (1967). Hitting time and hitting places for non-lattice recurrent random walks. *J. Math. Mech.* **17** 35–57. MR0215375

[24] PORT, S. C. and STONE, C. J. (1969). Potential theory of random walks on Abelian groups. *Acta Math.* **122** 19–114. MR0261706

[25] RACHEV, S. T. and SAMORODNITSKY, G. (1995). Limit laws for a stochastic process and random recursion arising in probabilistic modelling. *Adv. in Appl. Probab.* **27** 185–202. MR1315585

[26] SPITZER, F. (1964). *Principles of Random Walk*. D. Van Nostrand Co., Inc., Princeton, N.J.-Toronto-London. MR0171290

INSTITUTE OF MATHEMATICS
UNIVERSITY OF WROCŁAW
PL. GRUNWALDZKI 2/4
50-384 WROCŁAW
POLAND
E-MAIL: dbura@math.uni.wroc.pl